\documentclass{article}

\usepackage{amsmath} 
\usepackage{amssymb}
 \usepackage{amsthm} 
 \usepackage{enumerate} 
 \usepackage[mathscr]{eucal}
 \usepackage{tikz} 
\usetikzlibrary{arrows}
\usepackage{graphicx}
\usepackage{caption}
\usepackage{subcaption}
\usepackage{multicol}

\usepackage{fancyhdr}
\usepackage{epsfig}
\usepackage{latexsym}
\usepackage{microtype}
\usepackage{lastpage} 
\usepackage{hyperref}
\usepackage{multicol}
 \usepackage{tikz} 
 \usepackage{bbold}

\usepackage[all]{xypic}

\newcounter{cont}

\newtheorem{definition}{Definition}
\newtheorem{lemma}{Lemma}
\newtheorem{theorem}{Theorem}
\newtheorem{corollary}{Corollary}
\newtheorem{remark}{Remark}
\newtheorem{example}{Example}

\newcommand{\J}{\mathrm{J}}

\newcommand{\Stone}{\mathbf{B}\mathrm{ool}\mathbf{S}}
\newcommand{\StoneO}{\mathbf{B}\mathrm{ool}\mathbf{SO}_n}
\newcommand{\IStone}{\mathbf{S}\mathrm{t(I)}}
\newcommand{\Bool}{\mathbf{B}\mathrm{ool}\mathbf{A}}

\newcommand{\LMStone}{\mathbf{S}\mathrm{t(LM_{n+1})}}
\newcommand{\LMn}{\mathbf{LM}_{n+1}}
\newcommand{\BIn}{\mathbf{B}\mathrm{ool}\mathbf{I}_{n+1}}

\newcounter{cont1}

\title{ Mutually exclusive nuances of truth in Moisil logic} 
\author{Denisa Diaconescu and Ioana Leu\c{s}tean\\
{\small Department of Computer Science,} \\
{\small Faculty of Mathematics and Computer Science, University of Bucharest,}\\
{\small Academiei nr.14, sector 1, C.P. 010014,  Bucharest, Romania}\\   
{\small Emails: ddiaconescu@fmi.unibuc.ro,  ioana@fmi.unibuc.ro}}

\date{}

\begin{document}
\maketitle

\begin{abstract}

Moisil logic, having as algebraic counterpart \L ukasiewicz-Moisil algebras, provide an alternative  way to reason about vague information based on the following principle: a many-valued event is characterized by a family of Boolean events.   However,  using the original definition of \L ukasiewicz-Moisil algebra, the principle does not apply for subalgebras.   In this paper we identify an alternative and equivalent definition for the $n$-valued \L ukasiewicz-Moisil algebras, in which the determination principle is also saved for arbitrary subalgebras, which are characterized  by a Boolean algebra and a family of Boolean ideals.  As a consequence,  we prove a duality result  for the $n$-valued \L ukasiewicz-Moisil algebras, starting from the dual space of their Boolean center.  This leads us to a duality  for  MV$_n$-algebras, since  are equivalent to a subclass of $n$-valued \L ukasiewicz-Moisil algebras.

\end{abstract}

\section*{Introduction}

The first systems of many-valued logic are the $3$-valued and the $n$-valued  {\em \L ukasiewicz logic}   introduced by 
J. \L ukasiewicz  in the 1920's, while the infinite valued \L ukasiewicz logic was defined by J. \L ukasiewicz and A. Tarski in 1930   \cite{Lukasiewicz1920, Lukasiewicz-Tarski1930}. The investigation of the corresponding algebraic structures was a natural problem. The first who studied such an algebrization was Gr. C. Moisil who in 1941 introduced $3$ and $4$-valued {\em \L ukasiewicz algebras} \cite{Mosil1941}, and generalized later to the $n$-valued case \cite{Moisil1941} and the infinite case \cite{Moisil1972}. On the other hand, in 1958, C. C. Chang defined  {\em MV-algebras} \cite{Chang1958} as algebraic structures for \L ukasiewicz logic. As an example of A. Rose showed in 1965  that for $n \geq 5$ the \L ukasiewicz implication cannot be defined in an $n$-valued \L ukasiewicz algebra, the structures introduced by Moisil are not appropriate algebraic counterpart for \L ukasiewicz logic. In consequence, we are dealing with two different logical systems with different flavour: \L ukasiewicz logic, from one side, having MV-algebras as algebraic counterpart, and {\em Moisil logic}, from another side, having \L ukasiewicz algebras as corresponding algebras. Nowadays, we call \L ukasiewicz algebras by {\em \L ukasiewicz-Moisil algebras} and the standard monograph on these structures is \cite{BFGR}.

The proper subclass of \L ukasiewicz-Moisil algebras that correspond to $n$-valued \L ukasiewicz logic, i.e. {\em proper} \L ukasiewicz-Moisil algebras are characterized in \cite{Cignoli1982}. Since $MV_n$-algebras \cite{Grigolia1977} are the algebraic correspondent of the finite valued \L ukasiewicz logic, proper \L ukasiewicz-Moisil algebras and $MV_n$-algebras are categorical equivalent. The complex connections between \L ukasiewicz-Moisil algebras and MV-algebras are deeply investigated in \cite{IorgulescuI, IorgulescuII, IorgulescuIV}.

The main idea behind Moisil logic is that of {\em nuancing}: to a many-valued object we associate some Boolean objects, its {\em Boolean nuances}. We do not define an many-valued object by its Boolean nuances, but we characterize it through the nuances, we investigate its properties by reducing them to the study of some Boolean ones. Moisil logic is therefore derived from the classical logic by the idea of nuancing, mathematically expressed by a categorical adjunction. This is a general idea, that can be applied to any logical system, as pointed out in \cite{GG-Popescu2006}.

Moisil's {\em determination principle} plays a central role in the study of \L ukasiewicz-Moisil algebras and Moisil logic. At algebraic level it gives an efficient method for obtaining important results, lifting properties of Boolean algebras to the level of \L ukasiewicz-Moisil algebras (see Moisil's representation theorem, for example), while at logical level it gives an alternative way to reason about non-crisp objects, by evaluating some crisp ones.

The determination principle from the initial definition of \L ukasiewicz-Moisil algebras does not hold in general for subalgebras. The initial definition for \L ukasiewicz-Moisil algebras is given in terms of some lattice endomorphisms, called the {\em Chrysippian endomorphisms} by Moisil, i.e. the $\varphi$'s. Using  another family of unary  operations, i.e. the $\J$'s introduced in  \cite{Cignoli1982},   the determination principle for subalgebras can be proved  \cite{Ioana2008},   leading to the idea that  the Boolean nuances of a subalgebra are Boolean ideals. Moreover,  the alternative nuances $\J$'s  are {\em mutually exclusive}, or simply {\em disjoint}.

In this paper we introduce an  equivalent definition for \L ukasiewicz-Moisil algebras using the $\J$'s and we further investigate the properties of these operations.   We obtain a categorical equivalence that allow us to represent any \L ukasiewicz-Moisil algebra as a Boolean algebra and a finite family of Boolean ideals. As a consequnece,  we develop a duality for \L ukasiewicz-Moisil algebras starting from Boolean spaces and adding a family of open sets.  As a corollary, we obtain a duality for $MV_n$-algebras.

\medskip \medskip
The paper is organized as follows. In Section \ref{preliminaries} we recall the basic definitions and properties on Lukasiewicz-Moisil algebras. In particular in \ref{T(B)} we present the adjunction between \L ukasiewicz-Moisil algebras and Boolean algebra, the fundamental idea behind Moisil logic, while in \ref{Cignoli-duality} we recall the Stone-type duality  for \L ukasiewicz-Moisil algebras  \cite{cignoli1969}, which is developed starting from the dual space of a bounded distributive lattice. In Section \ref{alternative-definition} we introduce an alternative definition for \L ukasiewicz-Moisil algebras using the $\J$'s. In Subsection \ref{AdJ} we prove the fundamental logic adjunction theorem  via the J's, while in Subsection \ref{catEchiv} we prove a categorical equivalence for \L ukasiewicz-Moisil algebras. Sections \ref{sd1} and \ref{sd2} are devoted for Stone-type dualities for \L ukasiewicz-Moisil algebras and MV-algebras, respectively, starting from the simple dual space of Boolean algebras.

\section{\L ukasiewicz-Moisil algebras}\label{preliminaries}

In the sequel $n$ is a natural number and we use the notation $[n] := \{1,\ldots,n\}$.

\begin{definition}\label{definitionF}
A {\em \L ukasiewicz-Moisil algebra of order $n+1$} ($LM_{n+1}$-algebra, for short) is a structure of the form $$(L,\vee,\wedge,^*,\varphi_1,\ldots,\varphi_n,0,1)$$ such that $(L,\vee,\wedge,^*,0,1)$ is a De Morgan algebra, i.e. a bounded distributive lattice with a decreasing involution $^*$ satisfying the De Morgan property,  and $\varphi_1,\ldots,\varphi_n$\footnote{These operations are called the {\em Chrysippian endomorphisms.}} are unary operations on L such that the following hold:
\begin{list}
{(L\arabic{cont})}{\usecounter{cont}}
	\item $\varphi_i(x\vee y) = \varphi_i(x) \vee \varphi_i(y)$,
	\item $\varphi_i(x) \vee \varphi_i(x)^* = 1$,
	\item $\varphi_i \circ \varphi_j = \varphi_j$,
	\item $\varphi_i(x^*) = \varphi_{n+1-i}(x)^*$,
	\item if $i \leq j$ then $\varphi_i(x) \leq \varphi_j(x)$,
	\item if $\varphi_i(x) = \varphi_i(y)$, for all $i\in[n]$, then $x = y$,
\end{list}
for any $i,j\in[n]$ and $x,y\in L$.
\end{definition}

Property (L6) is called  {\em the determination principle} and the system (L1)-(L6) is equivalent to (L1)-(L5), (L7) and (L8), where
\begin{list}
{\em (L\arabic{cont})}{\usecounter{cont}\setcounter{cont}{6}\setlength{\leftmargin}{2cm}}
	\item $x \leq \varphi_n(x)$,
	\item $x \wedge \varphi_i(x)^* \wedge \varphi_{i+1}(y) \leq y$, for any $i\in[n-1]$.
\end{list}
Therefore, the class of $LM_{n+1}$-algebras is equational.

\begin{example}
The {\em canonical} $LM_{n+1}$-algebra is the structure $$(L_{n+1},\vee,\wedge,*,\varphi_1,\ldots,\varphi_n,0,1),$$ where $L_{n+1} := \{0,\frac{1}{n},\ldots,\frac{n-1}{n},1\}$, the lattice order is the natural one,
\begin{center}
$\frac{j}{n}^* := \frac{n-j}{n}$\hspace{0.5cm} and \hspace{0.5cm}$\varphi_i(\frac{j}{n}) :=  \left\{\begin{array}{rl}
								0, &\ \mathrm{if }\ i+j < n+1, \\
								1, &\ \mathrm{if }\ i+j \geq n+1\\
								\end{array}
			\right.,$
\end{center}
for any $0\leq j\leq n$ and $i\in [n]$.
\end{example}

%

The canonical $LM_2$-algebra has only one Chrysippian endomophism which, by the determination principle, it forced to be a bijection, making the canonical $LM_2$-algebra a Boolean algebra. Therefore the "overloaded" notation $L_2$ for Boolean algebras and $LM_2$-algebra is consistent.

By {\em Moisil's representation theorem}, 
any $LM_{n+1}$-algebra is isomorphic to a subdirect product of $LM_{n+1}$-subalgebras of $L_{n+1}$.

\begin{lemma}
In any $LM_{n+1}$-algebra $L$, the following hold, for any $x,y\in L$ and any $i,j\in [n]$:
\begin{multicols}{2}
\begin{list}
{(\arabic{cont})}{\usecounter{cont}}
	\item $\varphi_i(x\wedge y) = \varphi_i(x) \wedge \varphi_i(y)$,  
	\item $\varphi_i(x) \wedge \varphi_i(x)^* = 0$, 
	\item $\varphi_i(\varphi_j(x)^*) = \varphi_j(x)^*$, 
	\item $x \leq y$ iff $\varphi_i(x) \leq \varphi_j(x)$,
	\item $\varphi_1(x) \leq x \leq \varphi_n(x)$.
	\item[] $ $
\end{list}
\end{multicols}
\end{lemma}

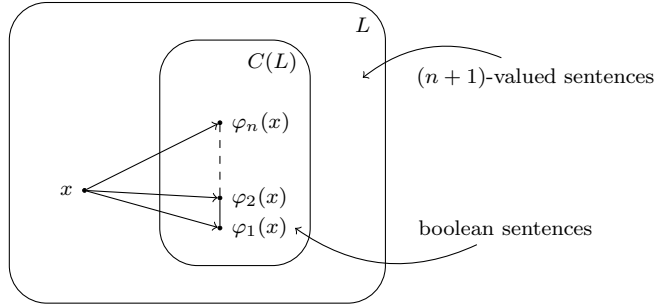
\begin{figure}[t]
\begin{center}
\begin{tikzpicture}
    [point/.style = {draw, circle,  fill = black, inner sep = .5pt}]		
    \draw [rounded corners=5mm] (1,1) rectangle (6,5);
    \draw [rounded corners=5mm] (3,1.5) rectangle (5,4.5);
    
    \node (C) at (4.5,4.2) {\footnotesize$C(L)$};
    \node (L) at (5.7,4.7) {\footnotesize$L$};
   
   \node (x) at (2,2.5) [point, label = {left:\footnotesize$x$}]{}; 
   \node (f1x) at (3.8,2) [point, label = {right:\footnotesize$\varphi_1(x)$}]{}; 
   \node (f2x) at (3.8,2.4) [point, label = {right:\footnotesize$\varphi_2(x)$}]{}; 
   \node (fnx) at (3.8,3.4) [point, label = {right:\footnotesize$\varphi_n(x)$}]{}; 
   
   \draw (f1x) -- (f2x);
   \draw[dashed] (f2x) -- (fnx);
   \draw[->] (x) -- (f1x);
   \draw[->] (x) -- (f2x);
   \draw[->] (x) -- (fnx);
   
   \node (n-valued) at (8,4) {\footnotesize $(n+1)$-valued sentences};
   \node (boolean) at (7.6,2) {\footnotesize boolean sentences};
   
   \path[->] (boolean) edge [bend left] (4.8,2);
   \path[->] (n-valued) edge [bend right] (5.7,4);
      
\end{tikzpicture}
\end{center}
\caption{The determination principle for $LM_{n+1}$-algebras.}
\label{determination principle}
\end{figure}

For each $LM_{n+1}$-algebra $L$, we define its {\em Boolean center} $C(L)$ as the set of all complemented elements of $L$, i.e. $C(L) = \{x\in L\ |\ x \vee x^*=1\}$. We can easily see that, for each $x\in L$, the following equivalences hold:
\begin{eqnarray*}
x \in C(L) & \mbox{iff} &  \mbox{ there exists }i\in [n]\mbox{ such that }\varphi_i(x) = x,\\
&\mbox{iff}& \mbox{ for all }i\in [n], \varphi_i(x) = x,\\
&\mbox{iff}& \mbox{ there exist }i\in [n]\mbox{ and }y\in L\mbox{ such that }\varphi_i(y) = x.
\end{eqnarray*}
Note that the determination principle can be represented as in Figure \ref{determination principle}.

\subsection{The Fundamental Logic Adjunction Theorem}\label{T(B)}
The logic having as algebraic semantics $LM_{n+1}$-algebras is called nowadays {\em Moisil logic}.
Moisil logic is derived from the classical logic by the {\em idea of nuancing}, mathematically expressed by a categorial adjunction. 

Let $(B,\vee,\wedge,^-,0,1)$ be a Boolean algebra. Let us consider the set
\begin{center}
$T(B) := \{(x_1,\ldots,x_n)\ |\ x_1 \leq \ldots \leq x_n\}$.
\end{center}
On the set $T(B)$ we can define an $LM_{n+1}$-algebra structure as follows: the lattice operations, as well as 0 and 1, are defined componentwise, and for each $(x_1,\ldots,x_n)\in T(B)$ and $i\in [n]$, we consider
	\begin{center}
	$(x_1,\ldots,x_n)^*:= (\overline{x_{n}},\ldots,\overline{x_{1}})$\hspace{.5cm} and \hspace{.5cm} $\varphi_i(x_1,\ldots,x_n) := (x_i,\ldots,x_i)$.
	\end{center}
Remark that the $LM_{n+1}$-algebras $L_{n+1}$ and $T(L_2)$ are isomorphic. 

Let $\LMn$ be the category of $(n+1)$-valued \L ukasiewicz-Moisil algebras and $\Bool$ be the category of Boolean algebras. Let 
\begin{center}
$C: \LMn \to \Bool$ \hspace{.5cm} and \hspace{.5cm} $T : \Bool \to \LMn$
\end{center}
 be two functors defined as follows:  for each $LM_{n+1}$-algebra, $C(L)$ is the Boolean center of $L$, and for each $LM_{n+1}$-morphism $f : L \to L'$, $C(f) : C(L) \to C(L')$ is the restriction and co-restriction of $f$ to $C(L)$ and $C(L')$; for each Boolean algebra $B$, $T(B)$ is the above $LM_{n+1}$-algebra, and for each Boolean morphism $g:B \to B'$, $T(g) : T(B) \to T(B')$ is defined by applying $g$ on each component of any $u\in T(B)$.
 
 
\begin{theorem}[The Fundamental Logic Adjunction Theorem]
The above functors $C$ and $T$ satisfy the following:
\begin{list}
{(\arabic{cont})}{\usecounter{cont}}
	\item $C$ is faithful and $T$ is full and faithful,
	\item $C$ is a left adjoint of $T$, where the unit and the counit are given by
	\begin{center}
	$\eta_L : L \to TC(L)$, $\eta_L(x)(i) = \varphi_i(x)$, for any $x\in L$ and any $i\in [n]$,
	
	$\epsilon_B : CT(B) \to B$, $\epsilon_B(u) = u(1)$, for all $u \in CT(B)$.
	\end{center}
	\item $\eta_L$ is an $LM_{n+1}$ embedding and $\epsilon_B$ is a Boolean isomorphism.
\end{list}
\end{theorem}

\subsection{Stone duality using the Chrysippian endomorphisms}\label{Cignoli-duality}

Two categories $\mathbf{C}$ and $\mathbf{D}$ are {\em dually equivalent}
if there exists a pair of contravariant functors $F \colon \mathbf{C} \to \mathbf{D}$ and $G \colon \mathbf{D} \to \mathbf{C}$
such that both $FG$ and $GF$ are naturally isomorphic with the corresponding identity functors, i.e.
for each object $C$ in $\mathbf{C}$ and $D$ in $\mathbf{D}$
there are isomorphisms $\eta_C \colon {GF}(C) \to C$ and $\kappa_D \colon {FG}(D) \to D$
such that
\begin{center}
$\xymatrix{
{GF}(C_1) \ar[d]_{\eta_{C_1}} \ar[r]^{{GF}(f)} & {GF}(C_2) \ar[d]^{\eta_{C_2}} \\
C_1 \ar[r]_f & C_2
}
\qquad
\xymatrix{
{FG}(D_1) \ar[d]_{\kappa_{D_1}} \ar[r]^{{FG}(g)} & {FG}(D_2) \ar[d]^{\kappa_{D_2}} \\
D_1 \ar[r]_g & D_2
}$
\end{center}
for each
$f \colon C_1 \to C_2$ in $\mathbf{C}$ and
$g \colon D_1 \to D_2$ in $\mathbf{D}$.

\medskip

A {\em Stone space} is a topological space $X$ such that: {\em (i)} $X$ is a compact $T_0$ space,
{\em (ii)} the set $K X$ of compact open subsets of $X$ is closed with respect to finite intersections and unions, and is a basis for the topology of $X$, and {\em (iii)} if $\mathcal{C} \subseteq K X$ is closed with respect to finite intersections and $F \subseteq X$ is a closed set such that $F \cap Y \neq \emptyset$, for every $Y \in \mathcal{C}$, then $F  \cap \bigcap_{Y\in \mathcal{C}} Y \neq \emptyset$. A map $f:X \to Y$ between two Stone spaces $X$ and $Y$ is called {\em strongly continuous} if $f^{-1}(A) \in K X$, for every $A \in K X$. We denote by $\Stone$ the category of Stone spaces and strongly continuous functions. The category of bounded distributive lattices is dually equivalent with the category $\Stone$.

A {\em Stone space with involution} is a couple $(X,h)$ such that: {\em (i)} $X$ is a Stone space, and {\em (ii)}
$h:X \to X$ is a function satisfying $h^2 = id_X$ and $X \setminus h(A) \in K X$, for every $A \in K X$.
We denote by  $\IStone$ the category of Stone spaces with involution, where the morphisms between $(X,g)$ and $(X',g')$ are the arrows $f:X \to X'$ in $\Stone$ such that $g' \circ f = f \circ g$. The category of De Morgan algebras  is dually equivalent with the category $\IStone$.

An {\em $(n+1)$-valued \L ukasiewicz-Moisil space} (an {\em $LM_{n+1}$-space}, for short) is a tuple $(X,h,h_1,\ldots,h_n)$ such that
\begin{list}
{\roman{cont})}{\usecounter{cont}}
	\item $(X,h)$ is a Stone space with involution,
	\item $\{h_i:X \to X\}_{i\in [n]}$ is a family of functions satisfying the conditions:
	\begin{list}
		{(\arabic{cont1})}{\usecounter{cont1}}
			\item $h_i$ is strongly continuous,
			\item $X \setminus h^{-1}_i(A) \in K X$, 
			\item $h_i \circ h_j = h_i$, 
			\item if $i\leq j$ then $h^{-1}_i(A) \subseteq h^{-1}_j(A)$, 
			\item if $h^{-1}_i(A) = h^{-1}_i(B)$, for all $i\in [n]$, then $A = B$, 
	\end{list}
	for any $i,j\in [n]$ and any  $A,B \in K X$,
	\item $h \circ h_i = h_{n+1-i}$ and $h_i \circ h = h_i$, for any $i\in [n]$,
\end{list}
A morphism between two $LM_{n+1}$-spaces $(Xh,h_1,\ldots,h_n)$ and $(X',h',h_1',\ldots,h_n')$ is an arrow $f:(X,h) \to (X',h')$ in $\IStone$ such that $f \circ g_i = g_i' \circ f$, for any $i\in [n]$. We denote by $\LMStone$ the category of $LM_{n+1}$-spaces.

\begin{theorem}\cite{cignoli1969}
 $\LMn$ and $\LMStone$ are dually equivalent.
\end{theorem}

\section{Mutually exclusive nuances of truth}\label{alternative-definition}

\subsection{Alternative definition}

\begin{definition}\label{definitionJ}
An $LM_{n+1}$-algebra is a structure of the form $$(L,\vee,\wedge,^*,\J_1,\ldots,\J_n,0,1)$$ such that $(L,\vee,\wedge,^*,0,1)$ is a De Morgan algebra and $\J_1,\ldots,\J_n$ are unary operations on L such that the following hold:
\begin{list}
{(J\arabic{cont})}{\usecounter{cont}}
	\item $\bigvee_{k = n-i+1}^n \J_k(x\vee y) = \bigvee_{k=n-i+1}^n(\J_k(x) \vee \J_k(y))$,
	\item $\J_i(x) \vee \J_i(x)^* = 1$,
	\item $\J_k (\J_i(x)) = 0$ and $\J_n(\J_i(x)) = \J_i(x)$,
	\item $\J_k(x^*) = \J_{n-k}(x)$ and $\J_n(x^*) = \bigwedge_{i=1}^n\J_i(x)^*$,
	\item $\J_l(x) \leq (\J_1(x) \vee \ldots \vee \J_{l-1}(x))^*$,
	\item if $\J_i(x) = \J_i(y)$, for all $i\in[n]$, then $x = y$,
\end{list}
for any $i,j\in[n]$, $k\in[n-1]$, $1 < l\leq n$ and $x,y\in L$.
\end{definition}

Note that for any $i,j\in [n]$ such that $i\neq j$, $\J_i(x)$ and $\J_j(x)$ are {\em mutually exclusive} (or {\em disjoint}, for short), i.e. $\J_i(x) \wedge \J_j(x) = 0$, for any $x\in L$. Indeed, assuming that $i<j$ and using (J5) and (J2), we have
\begin{eqnarray*}
\J_i(x) \wedge \J_j(x) = \J_i(x) \wedge \J_j(x) \wedge  \J_{j-1}(x)^* \wedge \ldots \wedge \J_1(x)^* = 0.
\end{eqnarray*}

\medskip
The two alternative definitions for $LM_{n+1}$-algebras are equivalent. Let $(L,\vee,\wedge,^*,\varphi_1,\ldots,\varphi_n,0,1)$ be a structure as in Definition \ref{definitionF}. As in \cite{Ioana2008}, we define the unary operations 
\begin{eqnarray}\label{defJ}
\J_n(x) := \varphi_1(x) \mbox{ and }\J_i(x) := \varphi_{n-i+1}(x) \wedge \varphi_{n-i}(x)^*, \mbox{ for }i \in [n-1], 
\end{eqnarray}
for any $x\in L$. Conditions (J2), (J3) and (J6) are already proved in \cite{Ioana2008}. Notice that, for any $i,k\in[n-1]$, we have
\begin{eqnarray*}
\bigvee_{j=i}^k\J_j(x) &=& [\varphi_{n-i+1}(x) \wedge \varphi_{n-i}(x)^*] \vee [\varphi_{n-i}(x) \wedge \varphi_{n-i-1}(x)^*] \vee \bigvee_{j=i+2}^k \J_j(x) \\
&=& [\varphi_{n-i+1}(x) \wedge \varphi_{n-i-1}(x)^*] \vee [\varphi_{n-i-1}(x) \wedge \varphi_{n-i-2}(x)^*] \vee \bigvee_{j=i+3}^k \J_j(x) \\
&=& [\varphi_{n-i+1}(x) \wedge (\varphi_{n-i+1}(x) \vee \varphi_{n-i-2}(x)^*) \wedge \varphi_{n-i-2}(x)^*]\vee \bigvee_{j=i+3}^k \J_j(x) \\
&=& [\varphi_{n-i+1}(x) \wedge  \varphi_{n-i-2}(x)^*]\vee \bigvee_{j=i+3}^k \J_j(x) = \ldots \\
&=& \varphi_{n-i+1}(x) \wedge \varphi_{n-k}(x)^*
\end{eqnarray*}
Therefore, for any $i,k\in[n-1]$, we also have
\begin{eqnarray*}
\bigwedge_{j=i}^k \J_j(x)^* = \varphi_{n-i+1}(x)^* \vee \varphi_{n-k}(x).
\end{eqnarray*}
Using the above equality, (J5) is obtained as follows, for any $1 < l \leq n$:
\begin{eqnarray*}
\J_l(x) \wedge \bigwedge_{j=1}^{l-1}\J_j(x)^* &=& \varphi_{n-l+1}(x) \wedge \varphi_{n-l}(x)^* \wedge (\varphi_n(x)^* \vee \varphi_{n-l+1}(x)) \\
&=& \varphi_{n-l+1}(x) \wedge \varphi_{n-l}(x)^* = \J_l(x)\\
\end{eqnarray*}
The first part of condition (J4) follows from \cite{Ioana2008}, while for the second part we have
\begin{eqnarray*}
\bigwedge_{j=1}^n\J_j(x)^* &=& \bigwedge_{j=1}^{n-1}\J_j(x)^* \wedge \J_n(x)^* 
= (\varphi_n(x)^* \vee \varphi_1(x)) \wedge \varphi_1(x)^* \\
&=&\varphi_n(x)^* \wedge \varphi_1(x)^* = \varphi_n(x)^* = \varphi_1(x^*) = \J_n(x^*) 
\end{eqnarray*}
Now let us show condition (J1) by induction over $i$. For $i=1$, the conclusion follows from (L1). Assume that (J1) holds for $i$ and let us prove it for $i+1$. First notice that, for any $2\leq i \leq n$, we have $\varphi_i(x)^*\vee \J_{n-i+1}(x) = \varphi_{i-1}(x)^*$. We have the following chain of equalities:
\begin{eqnarray*}
&& \bigvee_{k=n-(i+1)+1}^n \J_k(x \vee y) = \J_{n-i}(x \vee y) \vee \bigvee_{k = n-i+1}^n \J_k(x \vee y) \\
&=& [\J_{n-i}(x) \wedge \varphi_i(y)^*] \vee [\J_{n-i}(y)\wedge \varphi_i(x)^*]\vee \bigvee_{k=n-i+1}^n(\J_k(x) \vee \J_k(y)) \\
&=& [\J_{n-i}(x) \wedge \varphi_i(y)^*] \vee [\J_{n-i}(y)\wedge \varphi_i(x)^*] \vee [\J_{n-i+1}(x) \vee \J_{n-i+1}(y)] \vee\\
&& \bigvee_{k=n-i+2}^n(\J_k(x) \vee \J_k(y)) \\
&=& [(\J_{n-i}(x) \vee \J_{n-i+1}(y)) \wedge \varphi_{i-1}(y)^*] \vee [(\J_{n-i}(y) \vee \J_{n-i+1}(x)) \wedge \varphi_{i-1}(x)^*] \vee \\
&&\bigvee_{k=n-i+2}^n(\J_k(x) \vee \J_k(y)) \\
&=& \ldots \\
&=& [(\J_{n-i}(x) \vee \J_{n-i+1}(y) \vee \ldots \vee \J_{n-1}(y)) \wedge \varphi_1(y)^*] \vee \\
&&[(\J_{n-i}(y) \vee \J_{n-i+1}(x) \vee \ldots \vee \J_{n-1}(x)) \wedge\varphi_1(x)^*] \vee \J_n(x) \vee \J_n(y) \\
&=& [\J_{n-i}(x) \vee \ldots \vee \J_n(x)] \vee  [\J_{n-i}(y) \vee \ldots \vee \J_n(y)]  \\
&=& \bigvee_{k=n-(i+1)+1}^n(\J_k(x) \vee \J_k(y)).
\end{eqnarray*}

\medskip
\noindent Conversely, let $(L,\vee,\wedge,^*,\J_1,\ldots,\J_n,0,1)$ be a structure as in Definition \ref{definitionJ}. Again as in \cite{Ioana2008}, we define
\begin{eqnarray}\label{defF}
\varphi_i(x) := \bigvee_{k=n-i+1}^n\J_k(x),
\end{eqnarray}
for any $i\in [n]$ and $x\in L$. Condition (L5) follows immediately from the definition of $\varphi_i$'s, (L1) follows directly from (J1), while (L6) follows from (J6) noticing that $\J_i(x) = \varphi_{n-i+1}(x) \wedge \varphi_{n-i}(x)^*$, for any $i\in [n-1]$, and $\J_n(x) = \varphi_1(x)$. Using (J2), condition (L2) is obtained as follows:
\begin{eqnarray*}
\varphi_i(x) \vee \varphi_i(x)^* &=& \bigvee_{k=n-i+1}^n\J_k(x) \vee \bigwedge_{j=n-i+1}^n\J_j(x)^* \\
&=& \bigwedge_{j=n-i+1}^n (\bigvee_{k=n-i+1}^n \J_k(x) \vee \J_j(x)^*) = 1 \\
\end{eqnarray*}
Since for any $j\neq k$ and any $x\in L$, $\J_j(x)$ and $\J_k(x)$ are disjoint, it follows that $\J_j(x) \vee \J_k(x)^* = \J_k(x)^*$. Therefore, using (J2) and (J4), we get (L4):
\begin{eqnarray*}
\varphi_i(x^*) &=& \bigvee_{k=n-i+1}^n \J_k(x^*) = \bigvee_{k=n-i+1}^{n-1}\J_k(x^*) \vee \J_n(x^*) \\
&=&  \bigvee_{k=n-i+1}^{n-1} \J_{n-k}(x) \vee \J_n(x^*) = \bigvee_{j=1}^{i-1} \J_j(x) \vee \bigwedge_{k=1}^n\J_k(x)^* \\
&=& \bigwedge_{k=1}^n(\bigvee_{j=1}^{i-1}\J_j(x) \vee \J_k(x)^*) = \bigwedge_{k=i}^n(\bigvee_{j=1}^{i-1}\J_j(x) \vee \J_k(x)^*) \\
&=& \bigwedge_{k=i}^n\J_k(x)^* = (\bigvee_{k=i}^n\J_k(x))^* =\varphi_{n-i+1}(x)^*
\end{eqnarray*}
Finally, let us prove condition (L3) by induction over $j$. For $j=1$, using (J3) we have
\begin{eqnarray*}
\varphi_i(\varphi_1(x)) = \bigvee_{k=n-i+1}^n\J_k(\J_n(x)) = \J_n(x) = \varphi_1(x).
\end{eqnarray*}
Assume that (L3) holds for $j$ and let us prove it for $j+1$. Noticing that $\varphi_{j+1}(x) = \varphi_j(x) \vee \J_{n-j}(x)$ and using (J1) and (J3), we have
\begin{eqnarray*}
\varphi_i(\varphi_{j+1}(x)) &=& \bigvee_{k=n-i+1}^n \J_k(\varphi_{j+1}(x)) = \bigvee_{k=n-i+1}^n\J_k(\varphi_j(x) \vee \J_{n-j}(x)) \\
&=& \bigvee_{k=n-i+1}^n\J_k(\varphi_j(x)) \vee \bigvee_{k=n-i+1}^n\J_k(J_{n-j}(x)) \\
&=& \varphi_j(x) \vee \J_{n-j}(x) = \varphi_{j+1}(x).
\end{eqnarray*}

\subsection{The Fundamental Logic Adjunction Theorem  via J's}\label{AdJ}

Let $(B,\vee,\wedge,^-,0,1)$ be a Boolean algebra. For any family of elements $\{y_1,\ldots,y_n\}$ from $B$ we have
\begin{eqnarray*}
y_i \wedge y_j = 0,\ i\neq j & \mbox{ iff } & y_i\leq \overline{y_1} \wedge \ldots \wedge \overline{y_{i-1}},\ 1 < i \leq n,
\end{eqnarray*}
We call a family  of elements $\{y_1,\ldots,y_n\}$ with the above property {\em disjoint}. Let us define the set
\begin{center}
$J(B) := \{(y_1,\ldots,y_n)\in B^n \ |\ \{y_1,\ldots,y_n\}$ disjoint family of elements $\}$.
\end{center}

\begin{lemma}\label{bij} There exists a bijective correspondence between the sets $J(B)$ and $T(B)$.
\end{lemma}
\proof{
It is straightforward to check that the functions
\begin{eqnarray*}
f : J(B) \to T(B),& & f(y_1,\ldots,y_n) = (y_1,y_1\vee y_2,\ldots,y_1\vee \ldots \vee y_n), \\
g : T(B) \to J(B),& & g(x_1,\ldots,x_n) = (x_1,x_2\wedge \overline{x_1},\ldots,x_n\wedge \overline{x_{n-1}}) 
\end{eqnarray*}
define a bijective correspondence between the sets $J(B)$ and $T(B)$.
}

\medskip \medskip
As a consequence of Lemma \ref{bij}, $J(B)$ can be endowed with a structure of an $LM_{n+1}$-algebra $(J(B),\vee,\wedge,^*,\J_1,\ldots,\J_n,0,1)$. Consider the  $LM_{n+1}$-structure on $T(B)$ defined in Subsection \ref{T(B)}, $(T(B),\vee,\wedge,^*,\varphi_1,\ldots,\varphi_n,0,1)$. Using \ref{defJ}, the $\J$'s on the set $J(B)$ are defined as follows:
\begin{eqnarray*}
\J_i(y_1,\ldots,y_n) &:=& g(\varphi_{n-i+1}(f(y_1,\ldots,y_n)) \wedge \varphi_{n-i}(f(y_1,\ldots,y_n))^*) \\
\J_n(y_1,\ldots,y_n) &:=& g(\varphi_1(f(y_1,\ldots,y_n))) 
\end{eqnarray*}
for any $i\in[n-1]$. By simple computation, we obtain that
\begin{eqnarray*}
\J_i(y_1,\ldots,y_n) = (y_{n-i+1},0,\ldots,0) \mbox{ and } \J_n(y_1,\ldots,y_n) = (y_1,0,\ldots,0),
\end{eqnarray*} 
 for $i\in [n-1]$. In a similar way we get:
\begin{eqnarray*}
(y_1,\ldots,y_n)^* &= &(\bigwedge_{i=1}^n y^*_i, y_n^*,\ldots, y_2^*)\\
(y_1,\ldots,y_n)\vee (z_1,\ldots,z_n) &=& (w_1,\ldots,w_n)
\end{eqnarray*} 
where $w_1=y_1\vee z_1$,  $w_i=  (y_i\vee z_i)\wedge (y_{i-1}\vee z_{i-1})^*\wedge\cdots\wedge (y_1\vee z_1)^*$,  for  $i>1$.

\begin{remark}
The Fundamental Logic Adjunction Theorem from Subsection \ref{T(B)} can also be expressed  in terms of the functor $J$ defined in the obvious way.
\end{remark}

\subsection{A categorical equivalence for $LM_{n+1}$-algebras}\label{catEchiv}

We shall call a finite set $\{I_1,\ldots,I_{n-1}\}$ of ideals on a Boolean algebra $B$ with the property $I_i = I_{n-i}$, for any $i\in[n-1]$, an {\em $n$ symmetric sequence of Boolean ideals}. Consider the category $\BIn$ in which objects are tuples of the form $$(B,I_{n-1},\ldots,I_1),$$ where $B$ is a Boolean algebra and $\{I_1,\ldots,I_{n-1}\}$ is an $n$ symmetric sequence of Boolean ideals on $B$, and arrows are of the form $g:(B,I_{n-1},\ldots,I_1) \to (B',I_{n-1}',\ldots,I_1')$, where $g:B \to B'$ is a Boolean morphism and $g(I_i) \subseteq I_i'$, for any $i\in[n-1]$. 

Let us define the functor $$\Lambda : \LMn \to \BIn$$
as follows: for any $LM_{n+1}$-algebra $(L,\vee,\wedge,^*,\J_1,\ldots,\J_n,0,1)$, set $$\Lambda(L) = (C(L),J_{n-1}(L),\ldots,J_1(L)),$$ while for any $LM_{n+1}$-morphism $f:L\to L'$, set $\Lambda(f) : \Lambda(L) \to \Lambda(L')$ to be the co-restriction of $f$ to $C(L)$ and $C(L')$. The fact that $\Lambda(L)$ is an object in $\BIn$ follows from \cite[Proposition 5.2]{Ioana2008}.

Moreover, let us define the functor $$\Sigma : \BIn \to \LMn$$ as follows: for any object $(B,I_{n-1},\ldots,I_1)$ in $\BIn$, set 
\begin{center}
$\Sigma(B) = \{(y_1,\ldots,y_n) \in B^n \ |\ y_i\in I_{n-i+1}, \{y_1,\ldots,y_n\}$ disjoint elements$\},$
\end{center}
 and for each arrow  $g:(B,I_{n-1},\ldots,I_1) \to (B',I_{n-1}',\ldots,I_1')$ in $\BIn$, $\Sigma(g)$ is defined by applying $g$ on each component of any $u\in \Sigma(B)$. Applying \cite[Proposition 5.3]{Ioana2008} for the $LM_{n+1}$-algebra $J(B)$ defined in Subsection \ref{AdJ}, we obtain that $\Sigma(B)$ is an $LM_{n+1}$-algebra.

The functors $\Lambda$ and $\Sigma$ yield a categorical equivalence:

\begin{theorem}\label{th:cat}
The categories $\LMn$ and $\BIn$ are equivalent.
\end{theorem}

\section{Duality for $LM_{n+1}$-algebras using Boolean spaces}\label{sd1}

Let us recall the Stone duality for Boolean algebras. Let $\Bool$ be the category of Boolean algebras and Boolean homomorphisms and $\Stone$ the category of Boolean spaces (i.e. topological spaces that are Hausdorff and compact, and have a basis of clopen subsets) and continuous maps. The category $\Stone$ and $\Bool$ are dually equivalent via the functors
\begin{center}
$S^a:\Bool \to \Stone$ \hspace{.5cm} and \hspace{.5cm} $S^t : \Stone \to \Bool $.
\end{center}
The Boolean space $S^a(B)$ of a Boolean algebra $B$ is the topological space whose underlying set is the collection $X := Ult(B)$ of ultrafilters of $B$, and whose topology has a basis consisting of all sets of the form
$$N_b = \{U \in X\ |\ b \in U\},$$
for any $b\in B$. 
For every Boolean homomorphism $h:A \to B$, $S^a(f) : S^a(B) \to S^a(A)$ is defined as $S^a(f)(u) = h^{-1}(u)$, for every $u\in S^a(B)$. Conversely, if $X$ is a Boolean space, we consider $S^t(X) = co(X)$, the set of all clopen subsets of $X$, and for every continuous map $f:X \rightarrow Y$, $\varphi = S^t(f) : co(Y) \rightarrow co(X)$, $\varphi(N) = f^{-1}(N)$, for every $N\in co(Y)$. 

Note that for any $a,b\in B$, $N_{a\vee b} = N_a \cup N_b$, $N_{a\wedge b} = N_a \cap N_b$ and
\begin{eqnarray}
 \mbox{if } a\leq b \mbox{ then } N_a \subseteq N_b \nonumber
\end{eqnarray}
Via the Stone duality, for every ideal of a Boolean algebra we can associate an open set: if $I$ is an ideal of a Boolean algebra $B$, consider the open set in $S^a(B)$
\begin{eqnarray}
N_I  =  \bigcup\{N_a\ |\ a\in I\} =  \{U \in X\ |\ U \cap I \neq \emptyset\}. 
\end{eqnarray}
For every $b\in B$, we have $b\in I$ iff $N_b \subseteq N_I$. Conversely, for every open set of a Boolean space we can associate a Boolean ideal: if $O$ is an open subset of a Boolean space $X$, consider the ideal in $S^t(X)$
\begin{eqnarray}
I_O = \{b \in S^t(X)\ |\ N_b\subseteq O\}.
\end{eqnarray}

\medskip
In the following we provide a Stone-type duality for $LM_{n+1}$-algebras using Boolean spaces.
As we have seen in Subsection \ref{catEchiv} we can represent any $LM_{n+1}$-algebra $L$ as a Boolean algebra endowed with an $n$ symmetric sequence of Boolean ideals on it,  
$$(C(L),J_{n-1}(L),\ldots,J_1(L)),$$
the categories $\LMn$ and $\BIn$ being equivalent. Therefore we shall construct a Stone-type duality for the category $\BIn$.

\medskip
\begin{definition}
A {\em Boolean space with $n$ symmetric open sets} is a tuple $$(X,O_1,\ldots,O_{n-1}),$$ where  $X$ is a Boolean space and $O_1,\ldots,O_{n-1}$ are open sets in $X$ such that $O_i = O_{n-i}$, for any $i\in [n-1]$.
\end{definition}

Denote by $\StoneO$ the category of Boolean spaces with $n$ symmetric open sets with  continuous maps $f : (X,O_1,\ldots,O_{n-1}) \to (Y,U_1,\ldots,U_{n-1})$ such that $f^{-1}(U_i) \subseteq O_i$, for any $i\in [n-1]$.

\begin{theorem}\label{ThDual}
The categories $\BIn$ and $\StoneO$ are dually equivalent.
\end{theorem}
\proof{ We let define the functors
\begin{center}
$\Theta^a:\BIn \to \StoneO$ \hspace{.5cm} and \hspace{.5cm} $\Theta^t : \StoneO \to \BIn$. 
\end{center}
For every object $(B,I_{n-1},\ldots,I_1)$ in $\BIn$, define using (3),
\begin{center}
$\Theta^a(B,I_{n-1},\ldots,I_1) = (S^a(B),N_{I_1},\ldots,N_{I_{n-1}})$,
\end{center}
and for every arrow $g:(B,I_{n-1},\ldots,I_1) \to (B',I_{n-1}',\ldots,I_1')$ in $\BIn$, consider $\Theta^a(g) = S^a(g)$, {i.e. $\Theta^a(g)(u) = g^{-1}(u)$, for any $u\in S^a(B')$}. It is easy to check that $\Theta^a(g)^{-1}(N_{I_i}) \subseteq N'_{I'_{i}}$, for every $i\in[n-1]$, and therefore $\Theta^a(g)$ is an arrow in $\StoneO$. Conversely, for every Boolean space with $n$ symmetric open sets $(X,O_1,\ldots,O_{n-1})$, define using (4)
\begin{center}
$\Theta^t(X,O_1,\ldots,O_{n-1}) = (S^t(X), I_{O_{n-1}},\ldots,I_{O_1})$,
\end{center}
and for every arrow $f : (X,O_1,\ldots,O_{n-1}) \to (Y,U_1,\ldots,U_{n-1})$, define $\Theta^t(f) = S^t(f)$. Since $\Theta^t(f)(I_{U_i}) \subseteq I_{O_i}$, for any $i\in [n-1]$, $\Theta^t(f)$ is an arrow in $\BIn$.
}

\section{Duality for MV$_n$-algebras using Boolean spaces}\label{sd2}

An {\em MV-algebra} is an algebraic stucture $(A,\oplus, ^*, 0)$ of type $(2,1,0)$ such that $(A, \oplus, 0)$ is an abelian monoid $(x^*)^*=x$ and $x \oplus (y \oplus z)=(x\oplus y)\oplus z$. MV-algebras were defined by   \cite{Chang1958}  and they stay to \L ukasiewicz  $\infty$-valued logic as  Boolean algebras stay to classical logic. We refer to   \cite{cignoli-dottaviano-mundici} for an introduction in the theory of MV-algebras. 

If $A$ is an MV-algebra and $x,y\in A$ we set $x\odot y=(x^*\oplus y^*)^*$ and $1=0^*$. For any natural number $n$ we define
\begin{eqnarray*}
 0x=0, & x^0=1,\\
(n+1)x=(nx) \oplus x, & x^{n+1}=(x^n)\odot x.
\end{eqnarray*}

 The structures corresponding to the $n+1$-valued \L ukasiewicz logic were defined in \cite{Grigolia1977} under the name of $MV_n$-algebras  and they satisfy  the following additional properties for any $x\in A$  and $1 < j <n$  such that $j$ does not divide $n$:

$(n +1)x =nx$,

 $[(jx)\odot (x^* \oplus (( j −1)x)^*)]^n=0$.

In \cite{Cignoli1982} the {\em proper} $LM_{n+1}$-algebras are defined and they are those $LM_{n+1}$-algebras  adequate for the $n+1$-valued \L ukasiewicz logic. Moreover, $LM_{n+1}$-algebras are term equivalent with $MV_{n+1}$-algebras \cite{Iorgulescu1, IorgulescuII}. As in \cite{Ioana2008}, we call these structures
{\em Ł-proper}, in order to avoid the confusion with the usual terminology from universal algebra.  The \L -proper 
$LM_{n+1}$-algebras  are a full subcategory  of $\mathbf{LM}_{n+1}$.

\begin{remark}
In \cite[Section 3]{DiNola-Lettieri2000} the category $\mathbf{MV}_{n+1}$  is proved to be equivalent with a category $\mathbf{BM}_{n+1}$ whose objects are pairs $(B,R)$ such that $B$ is a Boolean algebra and $R\subset B^n$ satisfies some particular properties; the morphisms  form $(B_1,R_1)$ to $(B_2, R_2)$   in $\mathbf{BM}_{n+1}$  is a  morphism of Boolean algebras $f:B_1\to B_2$ such that $(x_1,\ldots, x_n)\in B_1$ implies $(f(x_1),\ldots, f(x_n))\in B_2$. Moreover, it is proved that  an object $(B,R)$ can  be characterized  by a sequence $I_1(R),\ldots, I_{n-1}(R)$  of Boolean ideals such that  
$J_i(L)\cap J_{i-k}(L)\subseteq J_{i}(L)$ for $2\leq i\leq n-2$ and $j<i$. In \cite{Ioana2008} this result is stated in the context of $LM_{n+1}$-algebras.
\end{remark}

\begin{lemma}\label{lem:help}
For an $LM_{n+1}$-algebra $L$ the following are equivalent:
\begin{enumerate}
\item[(a)] $L$ is \L -proper,
\item[(b)]  $J_i(L)\cap J_{k}(L)\subseteq J_{n-i+k-1}(L)$, for any $3\leq i\leq n-2$, $1\leq k\leq n-4$, $k<i$.
\end{enumerate}
\end{lemma}
\proof{ It  is straightforward by \cite[Proposition 5.11]{Ioana2008}.}

\medskip

Denote by $\mathbf{MV}_{n+1}$ the category of MV$_{n+1}$-algebras and by $\mathbf{BoolIMV}_{n+1}$  the full subcategory 
of $\BIn$ whose objects are  tuples of the form  
\begin{center}
$(B,I_{n-1},\ldots,I_1)$  such that  $I_i\cap I_{k}\subseteq I_{n-i+k-1}$,

 for any  $3\leq i\leq n-2$, $1\leq k\leq n-4$, $k<i$.
\end{center}
From Theorem \ref{th:cat} and  Lemma \ref{lem:help} we immediately infer the following result, which was proved directly  in \cite[Section 3]{DiNola-Lettieri2000}.

\begin{corollary}\label{CorMV}
The categories $\mathbf{MV}_{n+1}$ and  $\mathbf{BoolIMV}_{n+1}$ are equivalent.
\end{corollary}

The duality result for $MV_{n+1}$-algebras is now straightforward.  Denote by $\mathbf{BoolSOMV}_n$ the  full subcategory of category $\StoneO$  whose objects are Boolean spaces with $n$ symmetric open sets  
\begin{center}
$(X,O_1,\ldots,O_{n-1})$  such that  $O_i\cap O_{k}\subseteq O_{n-i+k-1}$,

 for any $3\leq i\leq n-2$, $1\leq k\leq n-4$, $k<i$.
\end{center}

\begin{theorem}

The categories $\mathbf{MV}_{n+1}$ and   $\mathbf{BoolSOMV}_n$  are dually equivalent.
\end{theorem}
\proof{ It is a direct consequence of Theorem \ref{ThDual} and Corollary \ref{CorMV}.

}

\newpage
\section{Conclusion}\label{conclusion}

Nuances of truth provide an alternative and robust way to reason about vague information: a many-valued object is uniquely determined by some Boolean objects, its nuances. However, a many-valued object cannot be recovered only from its Boolean nuances. This idea is mathematically expressed by a categorical adjunction between Boolean algebras and \L ukasiewicz-Moisil algebras.

Since the initial nuances of truth proposed by Moisil do not allow us to distinguish the  subalgebras, in this paper we explored a more expressible notion of nuances, namely mutually exclusive nuances of truth (or disjoint nuances of truth, for short).
 Apart from saving the determination principle for subalgebras, they  led us to a new Stone-type duality  for \L ukasiewicz-Moisil algebras.




\begin{thebibliography}{10}

\bibitem{BFGR}
V.~Boicescu, A.~Filipoiu, G.~Georgescu, and S.~Rudeanu.
\newblock {\em \L ukasiewicz-Moisil algebras}.
\newblock North-Holland, 1991.

\bibitem{Chang1958}
C.~C. Chang.
\newblock Algebraic analysis of many-valued logics.
\newblock {\em Trans. Amer. Math. Soc.}, 88:467--490, 1958.

\bibitem{cignoli1969}
R.~Cignoli.
\newblock {\em Algebras de Moisil de orden n}.
\newblock PhD thesis, Universidad Nacional de Sur, Bahia Blanca, 1969.

\bibitem{Cignoli1982}
R.~Cignoli.
\newblock Proper n-valued {\l}ukasiewicz algebras as s-algebras of
  {\l}ukasiewicz n-valued propositional calculi.
\newblock {\em Studia Logica}, 41:3--16, 1982.

\bibitem{cignoli-dottaviano-mundici}
R.~Cignoli, I.M.L. D'Ottaviano, and D.~Mundici.
\newblock {\em Algebraic Foundations of Many-Valued Reasoning}.
\newblock Kluwer Academic, 2000.

\bibitem{DiNola-Lettieri2000}
A.~Di~Nola and A.~Lettieri.
\newblock {One chain generated varieties of MV-algebras}.
\newblock {\em Journal of Algebra}, 225:667--697, 2000.

\bibitem{GG-Popescu2006}
G.~Georgescu and A.~Popescu.
\newblock {A common generalization for MV-algebras and {\L}ukasiewicz-Moisil
  algebras}.
\newblock {\em Archive for Mathematical Logic},
  45(8):947--981, 2006.




\bibitem{Grigolia1977}
R.S. Grigolia.
\newblock {Algebraic analysis of {\L}ukasiewicz-Tarski's logical systems}.
\newblock In R.~W\'{o}jcicki and G.~Malinowski, editors, {\em Selected Papers
  on {\L}ukasiewicz Sentensial Calculi}, pages 81--92. Osolineum, Wroclaw,
  1977.

\bibitem{IorgulescuI}
A.~Iorgulescu.
\newblock {Connections between MV$_n$ algebras and n-valued
  {\L}ukasiewicz-Moisil algebras - I}.
\newblock {\em Discrete Mathematics}, 181:155--177, 1998.

\bibitem{IorgulescuII}
A.~Iorgulescu.
\newblock {Connections between MV$_n$ algebras and n-valued
  {\L}ukasiewicz-Moisil algebras - II}.
\newblock {\em Discrete Mathematics}, 202:113--134, 1999.

\bibitem{IorgulescuIV}
A.~Iorgulescu.
\newblock {Connections between MV$_n$ algebras and n-valued
  {\L}ukasiewicz-Moisil algebras - IV}.
\newblock {\em The Journal of Universal Computer Science}, 6:139--154, 2000.

\bibitem{Ioana2008}
I.~Leu\c{s}tean.
\newblock {A determination principle for algebras of $n$-valued \L ukasiewicz
  logic}.
\newblock {\em Journal of Algebra}, 320:3694--3719, 2008.

\bibitem{Lukasiewicz1920}
J.~{\L}ukasiewicz.
\newblock O logice tr\'{o}jwarto\'{s}ciowej.
\newblock {\em Ruch. Filozoficzny}, 5(169-171), 1920.

\bibitem{Lukasiewicz-Tarski1930}
J.~{\L}ukasiewicz and A.~Tarski.
\newblock Untersuchungen \"{u}ber den aussagenkalk\"{u}l.
\newblock {\em C. R. Soc. Sc\'{e}ances Soc. Sci. Lettres Varsovie}, CL. III,
  23:30--50, 1930.

\bibitem{Mosil1941}
Gr.~C. Moisil.
\newblock Recherches sur les logiques non-chrysippiennes.
\newblock {\em Ann. Sci. Univ. Jassy}, 26:431--460, 1940.

\bibitem{Moisil1941}
Gr.~C. Moisil.
\newblock Notes sur les logiques non-chrysippiennes.
\newblock {\em Ann. Sci. Univ. Jassy}, 27(86-98), 1941.

\bibitem{Moisil1972}
Gr.~C. Moisil.
\newblock {\em Essais sur les logiques non-chrysippiennes}.
\newblock Editions de l'Academie de la Replublique Socialiste de Roumanie,
  Bucharest, 1972.

\end{thebibliography}

\newpage

\hypertarget{lastpage}{}
\lfoot{{\bf \copyright\hspace{0.01mm} Scientific Annals of Computer Science yyyy}}
\end{document}